\input amstex
\documentstyle{amsppt}
\magnification=\magstep1 \NoRunningHeads

\topmatter

\title  Spectral multiplicities for infinite~measure preserving transformations
 \endtitle
\author  Alexandre~I.~Danilenko and Valery V. Ryzhikov
\endauthor

\thanks
The second named author was supported in part by the Programme of Support of Leading Scientific Schools of RF (grant no. 3038.2008.1).
\endthanks

\abstract
Each subset $E\subset\Bbb N$ is realized as the set of essential values of the multiplicity function for  the Koopman operator of an ergodic conservative infinite measure preserving transformation.
\endabstract

\address
 Institute for Low Temperature Physics
\& Engineering of National Academy of Sciences of Ukraine, 47 Lenin Ave.,
 Kharkov, 61164, UKRAINE
\endaddress
\email alexandre.danilenko\@gmail.com
\endemail

\address
 Department of Mechanics and Mathematics, Lomonosov
Moscow State University, GSP-1, Leninskie Gory, Moscow, 119991, Russian
Federation
\endaddress
\email vryzh\@mail.ru
\endemail

\NoBlackBoxes
\endtopmatter
\document

\head 0. Introduction
\endhead
Let $T$ be an ergodic conservative invertible measure preserving transformation  of a
$\sigma$-finite standard measure space $(X,\goth B,\mu)$.
Consider an associated unitary (Koopman) operator $U_T$ in the Hilbert space $L^2(X,\mu)$:
$$
U_Tf:=f\circ T.
$$
In the case of finite measure $\mu$, the operator $U_T$ is usually considered
only in the orthocomplement to the subspace of constant functions.
A general question of the spectral theory of dynamical systems is
$$
\text{{\it to find out which unitary operators can be realized as Koopman operators.}}\tag0-1
$$
Several
 particular cases of \thetag{0-1} are well known in the case of finite $\mu$: Banach problem on simple
Lebesgue spectrum, Kolmogorov problem on group property of spectrum
\cite{KVe}, \cite{St}, Rokhlin problem on homogeneous spectrum and, more
generally, spectral multiplicity problem.
We state the latter one as follows.
Denote by $\Cal M(T)$ the set of essential values for the spectral
multiplicity function of $U_T$. Then
$$
\text{\it which subsets of $\Bbb N$ are realizable as $\Cal M(T)$ for an ergodic $T$?}\tag0-2
$$
Despite a significant  progress achieved in works of many authors \cite{Os}, \cite{Ro1}, \cite{Ro2}, \cite{G--L}, \cite{KwL}, \cite{Ka}, \cite{Ag1}, \cite{Ry1}, \cite{Ag2}, \cite{Da2},  \cite{Ry2}, \cite{Ag3}, \cite{KaL}, \cite{Da4}, \cite{Ry3} etc., this long-standing basic question of the spectral theory of finite measure preserving dynamical systems remains open.

In the present paper we  consider \thetag{0-2} in the class of {\it infinite} measure preserving conservative ergodic transformations.
It turns out that \thetag{0-2} can be solved completely in this class.

\proclaim{Theorem 0.1}
Given any subset $E\subset\Bbb N$, there is an ergodic multiply recurrent infinite measure preserving transformation $T$ such that
$\Cal M(T)=E$.
\endproclaim

If  $E\ni p+1$ for some positive integer $p$, we seek  the desired example in the form $T^{\odot p}\times T_\alpha$, where $T^{\odot p}$ is a symmetric power of an appropriate transformation $T$ and $T_\alpha$ is a certain compact extension of $T$.
If $E=\{p+1\}$ then the desired transformation is $T^{\odot p}\times T$.
To implement  this idea we adapt the techniques developed in a recent paper \cite{Da4} (which, in turn, absorbed the techniques from nearly all aforementioned papers) for probability preserving systems to the infinite setup.

We note that in the infinite measure preserving case we encounter some specific  phenomena which are absent in the probability preserving case.
For instance, ergodicity is not a spectral property.
Besides,  there is no a {\it good} definition for weak mixing.
Indeed, for each $p>0$, there is a transformation $T$ such that the $p$-fold Cartesian power $T^{\times p}$ is ergodic but $T^{\times (p+1)}$ is not \cite{KPa}.
Next, we recall that a transformation $T$ on $(X,\mu)$ is called {\it multiply recurrent} if for each subset $A\subset X$ of positive measure and each $p>0$, there exists $k>0$ such that $\mu(A\cap T^kA\cap\cdots \cap T^{kp}A)>0$. This concept  is a natural strengthening  of conservativeness (that corresponds to $p=1$).
Furstenberg showed that if $\mu(X)<\infty$ then each $\mu$-preserving transformation is multiply recurrent \cite{Fu}. However if $\mu(X)=\infty$ then there are ergodic transformations (even with all Cartesian powers ergodic) which are not multiply recurrent \cite{AFS}. For more information on these and many other {\it infinite} counterexamples we refer to \cite{Aa}, \cite{Da3}, \cite{DaS2} and references therein.

The reason why the proof of Main Theorem does not work  in the case of probability measure is that the factors $T^{\odot p}$ and $T_\alpha$ of
$T^{\odot p}\times T_\alpha$  have (in the probability case) one-dimensional invariant subspaces of constants which add ``superfluous'' multiplicities.

\head 1.  Rokhlin problem on multiplicities~for infinite measure preserving maps
\endhead

Rokhlin problem on multiplicities can be stated as follows
\roster
\item"---"
{\it given $n>1$, is there an ergodic transformation with homogeneous spectrum of multiplicity $n$?}
\endroster
This particular case of \thetag{0-1} plays an important role in the proof
of Main Theorem. We note that in the finite measure preserving case Rokhlin
problem was solved in \cite{Ag1} and \cite{Ry1} for $n=2$.
 For an arbitrary $n$ it was solved in \cite{Ag2}
and in a constructive way in \cite{Da2}.
To solve Rokhlin problem in the
infinite measure preserving case it is enough to consider natural factors
of Cartesian powers (cf. \cite{Ka}, \cite{Ag3}).

\proclaim{Lemma 1.1}
Let $\Gamma$ be a subgroup in $\goth S_k$.
Let $V$ be a unitary with a simple continuous spectrum such that $V^{\odot k}$ has a simple spectrum.
Denote by $V^{\otimes k}_\Gamma$ the restriction  of $V^{\otimes k}$ to the subspace of $\Gamma$-invariant tensors. Then
$V^{\otimes k}_\Gamma$ has a homogeneous spectrum of multiplicity
$k!/\# \Gamma$.
\endproclaim
\demo{Proof}
Let $\sigma$ stand for  a measure of maximal spectral type of $V$.
 Consider a homomorphism $\pi:\Bbb T^k\ni(z_1,\dots,z_k)\mapsto z_1\cdots z_k\in \Bbb T$.
Since $\sigma$ is continuous,
$\goth S_k$ acts freely on almost
every fiber $\pi^{-1}(z)$ furnished with the conditional measure $\sigma^k\restriction \pi^{-1}(z)$.
The operator  $V^{\odot k}$ has a simple spectrum.
Hence almost all of the conditional measures are concentrated at $k!$ points (the collection of these points depends on fiber) of positive measure.
The action of $\Gamma$ on a fiber partitions this set into $k!/\Gamma$ orbits of equal cardinality.
Thus the operator $V^{\otimes k}_\Gamma$ is unitarily equivalent to a direct sum of $k!/\Gamma$ copies of $V^{\odot k}$ each of which has a simple spectrum.
 \qed
\enddemo

\proclaim{Corollary 1.2} Let $T$ be an ergodic conservative infinite measure preserving conservative map such that $U_T^{\odot k}$ has a simple spectrum.
Then $\Cal M({T^{\odot (k-1)}\times T})=\{k\}$.
\endproclaim
\demo{Proof}
Since $T$ is ergodic of type $II_\infty$, it follows that $U_T$ has a continuous spectrum.
It remains to note that $U_{T^{\odot (k-1)}\times T}=U_T^{\odot (k-1)}\otimes U_T=(U_T)^{\otimes k}_{\goth S_{k-1}}$ and apply Lemma~1.1.
\qed
\enddemo

We note however that Corollary~1.2 does not solve the Rokhlin problem since it is unclear whether the  transformation $T^{\odot (k-1)}\times T$ is ergodic or not.
A complete solution will appear in the proof of Theorem~0.1.

\head 2. Preliminaries and notation
\endhead
Given a countable subset $A$, we denote by $\# A$ the cardinality of $A$.
Let $G$ be a countable Abelian group, $H$ a subgroup of $G$ and $v:G\to G$ a group automorphism.
As usual, $\widehat G$ and $\widehat v$ denote the dual group
and the dual automorphism  respectively.
We set
$$
\align
M_v^h&:= \#(\{v^i(h)\mid i\in\Bbb Z\}\cap H),\\
L(G,H,v)&:=\{M_v^h\mid h\in H\setminus\{0\}\},\\
\Cal G &:=\{a\in\widehat G\mid\exists\text{$p>0$ with } \widehat v^p(a)=a\},\\
l_g(a) &:=\frac 1p\sum_{i=0}^{p-1}a(v^i(g))\quad\text{for all $a\in\Cal G$ with $\widehat v^p(a)=a$ and $g\in G$}.
\endalign
$$
We  state in this section two auxiliary lemmata which will play the key role when we  compute the spectral multiplicities of the transformations under consideration.

\proclaim{Lemma 2.1} Given any subset $E\subset\Bbb N$, there exist
a countable Abelian group $G$, a subgroup $H\subset G$ and an automorphism $v:G\to G$ such that $E=L(G,H,v)$. Moreover, the following properties are satisfied:
\roster
\item"\rom{(i)}"
the subgroup $\Cal G$ is countable and dense in $\widehat G$,
\item"\rom{(ii)}"
if $g_1,g_2\in G$ and $v^i(g_1)\ne g_2$ for all $i\in\Bbb Z$ then there is $a\in\Cal G$  such that $l_{g_1}(a)\ne l_{g_2}(a)$.
\endroster
\endproclaim

For the proof of Lemma~2.1 we refer to  \cite{Da4}.
The following lemma  was proved in \cite{KaL} under  slightly stricter conditions.
We give here an alternative short proof of it.

\proclaim{Lemma 2.2} Let $V$ and $W$ be unitary operators with simple spectrum in Hilbert spaces $\Cal H$ and $\widetilde{\Cal H}$ respectively. Assume  that for each $i=1,\dots,k$, there are two sequences $n_t^{(i)}\to\infty$, $m_t^{(i)}\to\infty$ and three complex numbers
 $\widetilde\kappa_i\ne\kappa_i$ and $\delta_i$ such that the following weak limits exist:
\roster
\item"\rom{(i)}"
$V^{n_t^{(i)}}\to  \kappa_iI+\delta_iV^*$, $W^{n_t^{(i)}}\to  \kappa_iI+\delta_iW^*$,
\item"\rom{(ii)}"
$V^{m_t^{(i)}}\to \kappa_iI+\delta_i V^*$,
$W^{m_t^{(i)}}\to \widetilde\kappa_iI+\delta_i W^*$.
\endroster
 Assume, moreover, that $\#\{\kappa_1/\delta_1,\dots,\kappa_k/\delta_k\}=k$.
Then $V^{\odot k}\otimes W$ has a simple spectrum.
\endproclaim

\demo{Proof}
Without loss of generality we can think that $\Cal H=L^2(\Bbb T,\sigma_V)$,
$\widetilde{\Cal H}=L^2(\Bbb T,\sigma_W)$ and $Vf(z)=zf(z)$, $Wg(z)=zg(z)$ for all $f\in\Cal H$, $g\in\widetilde{\Cal H}$, $z\in\Bbb T$.
Here $\sigma_V$ and $\sigma_W$ stand for  measures of maximal spectral type of $V$ and $W$ respectively.
Then
$$
\align
&\Cal H^{\odot k}\otimes\widetilde{\Cal H}=L^2_{\text{sym}}(\Bbb T^k,\sigma_V^k)\otimes L^2(\Bbb T,\sigma_W)\quad\text{and}\\
&(V^{\odot k}\otimes W)f(z_1,\dots,z_k,z)=z_1\cdots z_kzf(z_1,\dots,z_k,z)
\endalign
$$
for all $f\in \Cal H^{\odot k}\otimes\widetilde{\Cal H}$.
Let $\Cal Z$  stand for the $(V^{\odot k}\otimes W)$-cyclic space generated by the constant function $1\in \Cal H^{\odot k}\otimes\widetilde{\Cal H}$.

Denote by $\Cal A$ the von Neumann algebra  generated by $V^{\odot k}\otimes W$. We consider  elements of $\Cal A$ as  bounded functions on $\Bbb T^{k+1}$ which are invariant under any permutation of the  $k$ first coordinates.
From (i) and (ii) we deduce that the following two functions
$$
\align
(z_1,\dots,z_k,z) &\mapsto (\kappa_i+\delta_i z)\cdot\prod_{l=1}^k(\kappa_i +\delta_i z_l)\\
(z_1,\dots,z_k,z) &\mapsto (\widetilde\kappa_i+\delta_i z)\cdot\prod_{l=1}^k(\kappa_i +\delta_iz_l)
\endalign
$$
are in $\Cal A$ for each $i=1,\dots,k$.
Since $\kappa_i\ne\widetilde\kappa_i$, it follows that the function
$$
(z_1,\dots,z_k,z) \mapsto \prod_{l=1}^k\bigg(\frac{\kappa_i}{\delta_i} +z_l\bigg)=
\sum_{l=0}^k
\bigg(\frac{\kappa_i}{\delta_i}\bigg)^l
P_l(z_1,\dots,z_k)
$$
is in $\Cal A$ for each $i=1,\dots,k$.
Hence $P_0,\dots,P_k$ are all in $\Cal A$ (by the property of Vandermond determinant and a condition of the lemma).
It is easy to see that the polynomials $P_0,\dots,P_k$ generate the entire algebra $\Cal P_{\text{sym}}(k)$ of symmetric polynomials in $k$ variables.
Since $\Cal Z$
 is invariant under $\Cal A$ and the linear subspace
$\Cal P_{\text{sym}}(k)1$ is dense in $L^2_{\text{sym}}(\Bbb T^k,\sigma_V^k)$, we then obtain that
$
\Cal Z\supset L^2_{\text{sym}}(\Bbb T^k,\sigma_V^k)\otimes 1.
$
Hence
$$
\Cal Z\supset L^2_{\text{sym}}(\Bbb T^k,\sigma_V^k)\otimes z^n
$$
for all $n\in\Bbb Z$.
This yields
$
\Cal Z=L^2_{\text{sym}}(\Bbb T^k,\sigma_V^k)\otimes L^2(\Bbb T,\sigma_W).
$
\qed
\enddemo

Let $E$ be any subset of $\Bbb N$.
Passing to the dual objects in the statement of Lemma~2.1 we obtain that there exist a compact Polish Abelian group $K$,  a closed subgroup $H$ of $K$ and  a continuous automorphism $v$ of $K$ such that
$$
E=L(\widehat K,\widehat{K/H}, \widehat v)
 $$
 and the conditions (i)--(iii) of Lemma~2.1 are satisfied.
The subset of $v$-periodic points in $K$ will be denoted by $\Cal K$.
We keep this notation till the end of the paper.

Let $T$ be an ergodic transformation of $(X,\mu)$. Denote by $\Cal R\subset X\times X$ the $T$-orbit equivalence relation.
 A measurable map $\alpha$ from $\Cal R$ to  $K$ is called a {\it cocycle} of $\Cal R$ if
$$
\alpha(x,y)\alpha(y,z)=\alpha(x,z)\quad \text{for all }(x,y), (y,z)\in\Cal R.
$$
 Two cocycles $\alpha,\beta:\Cal R\to K$ are {\it cohomologous} if
there are a measurable map $\phi:X\to K$ and a $\mu$-conull subset $B\subset X$ such that
$$
\alpha(x,y)=\phi(x)\beta(x,y)\phi(y)^{-1}\quad\text{for all }(x,y)\in \Cal R\cap (B\times B).
$$
If a transformation $S$ commutes with $T$ (i.e. $S\in C(T)$) then a cocycle $\alpha\circ S:\Cal R\to K$ is well defined by $\alpha\circ S(x,y):=\alpha(Sx,Sy)$.

Let $\lambda_{K/H}$ stand for Haar measure on ${K/H}$.
We denote by $T_{\alpha,H}$ the following transformation of the space
$(X\times K/H,\lambda_{K/H})$:
$$
T_{\alpha,H}(x,k+H):=(Tx,\alpha(Tx,x)+k+H).
$$
It is called a {\it skew product} extension of $T$.
For brevity, $T_{\alpha,\{0\}}$ will be denoted by  $T_\alpha$.
As in the case of finite measure preserving transformations we have a decomposition of  $U_{T_{\alpha,H}}$ into  orthogonal sum $U_{T_{\alpha,H}}=\bigoplus_{\chi\in\widehat{K/H}}U_{T,\chi}$, where  $U_{T,\chi}$ is the following unitary in $L^2(X,\mu)$:
$$
U_{T,\chi}f(x)=\chi(\alpha(Tx,x))f(Tx), \quad x\in X.
$$

\remark{Remark \rom{2.3}}
It is straightforward to verify that  if
$$
\alpha\circ S\text{ \ is cohomologous to \ }  v\circ\alpha
$$
for some $S\in C(T)$ then $U_{T,\chi}$ and $U_{T,\widehat v^i(\chi)}$ are unitarily equivalent for each $i\in\Bbb Z$.
\endremark

\head 3. Construction  of the base transformation and the cocycle
\endhead

To prove  Theorem~0.1 we will need rank-one transformations and their compact extensions.
 It is well known that the rank-one transformations have a simple spectrum.
We will construct them using  the $(C,F)$-construction (see \cite{dJ}, \cite{Da1}--\cite{Da3}).
We now briefly outline its formalism.
Let two sequences $(C_n)_{n>0}$ and $(F_n)_{n\ge 0}$ of finite subsets in $\Bbb Z$ are given such that:
\roster
\item"---"
$F_n=\{0,1,\dots,h_n-1\}$, $h_0=1$, $\# C_n>1$, $0\in C_n$,
\item"---" $F_n+C_{n+1}\subset F_{n+1}$,
\item"---" $(F_{n}+c)\cap (F_n+c')=\emptyset$ if $c\ne c'$, $c,c'\in C_{n+1}$,
\item"---" $\lim_{n\to\infty}\frac{h_n}{\#C_1\cdots\# C_n}=\infty$.
\endroster
Let $X_n:=F_n\times C_{n+1}\times C_{n+2}\times\cdots$. Endow this set with the (compact Polish) product topology. The  following map
$$
(f_n,c_{n+1},c_{n+2})\mapsto(f_n+c_{n+1},c_{n+2},\dots)
$$
is a topological embedding of $X_n$ into $X_{n+1}$. We now set $X:=\bigcup_{n\ge 0} X_n$ and endow it with the (locally compact Polish) inductive limit topology. Given $A\subset F_n$, we denote by $[A]_n$ the following cylinder: $\{x=(f,c_{n+1},\dots,)\in X_n\mid f\in A\}$. Then $\{[A]_n\mid A\subset F_n, n>0\}$ is the family of all compact open subsets in $X$.
It forms a base of the topology on $X$.
Denote by $\goth B$ the Borel $\sigma$-algebra generated by this topology.

Let  $\Cal R$ stand for the {\it tail} equivalence relation on $X$: two points $x,x'\in X$ are $\Cal R$-equivalent if there is $n>0$ such that $x=(f_n,c_{n+1},\dots),\ x'=(f_n',c_{n+1}',\dots)\in X_n$ and $c_m=c_m'$ for all $m>n$.
Recall that a measure  on $(X,\goth B)$ is  {\it $\Cal R$-invariant} if it is invariant under each measurable transformation $S:X\to X$ whose graph is contained in $\Cal R$.
An $\Cal R$-invariant measure $\mu$ is called $\Cal R$-{\it ergodic} if every $\Cal R$-saturated measurable subset of $X$ is either $\mu$-null or $\mu$-conull.
There is only one   non-atomic infinite $\sigma$-finite measure $\mu$ on $X$ which is invariant (and ergodic) under $\Cal R$ and such that $\mu(X_0)=1$.

Now we define a transformation $T$ of $(X,\mu)$ by setting
$$
T(f_n,c_{n+1},\dots):=(1+f_n,c_{n+1},\dots )\text{ whenever }f_n<h_n-1,\ n>0.
$$
This formula defines $T$ partly on $X_n$. When $n\to\infty$, $T$ extends to the entire $X$ (minus countably many points) as a $\mu$-preserving transformation. Moreover, the $T$-orbit equivalence relation  coincides  with $\Cal R$ (on the subset where $T$ is defined). We call $T$ {\it the $(C,F)$-transformation} associated with $(C_{n+1},F_n)_{n\ge 0}$.

We recall a concept of a $(C,F)$-cocycle (see \cite{Da2}).
 Given a sequence of maps $\alpha_n:C_n\to K$, $n=1,2,\dots$, we first define a Borel cocycle $\alpha:\Cal R\cap (X_0\times X_0)\to K$ by setting
$$
\alpha(x,x'):=\sum_{n>0}(\alpha_n(c_n)-\alpha_n(c_n')),
$$
whenever $x=(0,c_1,c_2,\dots)\in X_0$,  $x'=(0,c_1',c_2',\dots)\in X_0$ and $(x,x')\in\Cal R$. To extend $\alpha$ to the entire $\Cal R$, we first define a map $\pi:X\to X_0$ as follows. Given $x\in X$, let $n$ be the least positive integer such that $x\in X_n$. Then $x=(f_n,c_{n+1},\dots)\in X_n$. We  set
$$
\pi(x):=(\underbrace{0,\dots,0}_{n+1\text{ times}}, c_{n+1}, c_{n+2},\dots)\in X_0.
$$
Of course, $(x,\pi(x))\in\Cal R$. Now for each pair $(x,y)\in\Cal R$, we let
$$
\alpha(x,y):=\alpha(\pi(x),\pi(y)).
$$
It is easy to verify that $\alpha$ is a well defined cocycle of $\Cal R$ with values in $K$. We call it {\it the $(C,F)$-cocycle associated with} $(\alpha_n)_{n=1}^\infty$.

Let $\bar z=(z_n)_{n=1}^\infty$ be a sequence of integers.
The following statement is an {\it infinite} analogue of \cite{Da2, Lemma~4.11}.
A ``spectral meaning'' of it was explained in~Remark~2.3.

\proclaim{Lemma 3.1}  Suppose that
$
\sum_{n=1}^\infty\frac{\#(C_{n}\triangle(C_{n}-z_{n}))}{\# C_{n}}<\infty.
$
We  set
 $$
X^{\bar z}_n:=\{0,1,\dots, h_n-z_1-\cdots-z_n\}\times\prod_{m>n}(C_{m}\cap(C_m-z_m))\subset X_n.
$$
Then a transformation $S_{\bar z}$ of $(X,\mu)$ is well defined by setting
$$
S_{\bar z}(x):=(z_1+\cdots+z_n+f_n, z_{n+1}+c_{n+1}, z_{n+2}+c_{n+2},\dots)\tag3-1
$$
for all $x=(f_n,c_{n+1},c_{n+2},\dots)\in X^{\bar z}_n$, $n=1,2,\dots$.
Moreover, $S_{\bar z}$ commutes with $T$ and $U_T^{z_1+\cdots+z_n}\to U_{S_{\bar z}}$ weakly as $n\to\infty$.

Let $\alpha$ be the $(C,F)$-cocycle associated with a sequence of maps
$\alpha_n:C_n\to K$, $n=1,2,\dots$.
Let $C_n^\circ:=\{c\in C_n\cap (C_n-z_n)\mid \alpha_n(c+z_n)=v(\alpha_n(c))\}$.
If
$$
\sum_{n>0}(1-\# C_n^\circ/\#C_{n})<\infty\tag3-2
$$
then the cocycle $\alpha\circ S_{\bar z}$ is cohomologous to $v\circ\alpha$.
\endproclaim

\demo{Proof}
We need to verify that the subset $D:=\bigcup_{n>1}X_n^{\bar z}$ is of full measure in $X$.
Fix $n>0$.
We first note that
$$
(X_n,\mu\restriction X_n)=(F_n,\nu_n)\otimes\bigotimes_{m>n}(C_m,\kappa_m),
$$
where $\kappa_m$ is the equidistributed probability measure on $C_n$
and $\nu_n$ is an equidistributed finite measure on $F_n$ \cite{Da1}.
It follows from Borel-Cantelli lemma that
$$
X_n\setminus D\subset\bigcup_{m>n}[\{h_m-z_1-\cdots-z_m+1,\dots,h_m-1\}]_m
$$
up to a subset of $\mu$-null measure.
Therefore
$$
\mu(X_n\setminus D)\le\sum_{m>n}\frac{z_1+\cdots+z_m}{h_m}\mu(X_m)=
\sum_{m>n}\frac{z_1+\cdots+z_m}{\#C_1\cdots\# C_m}\mu(X_0)<2\sum_{m>n}\frac{z_m}{\# C_m}.
$$
Since $\sum_{m>1}{z_m}/{\# C_m}<\infty$, it follows that $\mu(X_n\setminus D)\to 0$ as $n\to\infty$.

The second claim can be shown by an obvious modification of the proof of Lemma~4.11 from \cite{Da2}. \qed
\enddemo

We will construct some special $(C,F)$-transformation and its cocycle with values in $K$. Fix a partition
$$
\Bbb N=\bigsqcup_{a\in \Cal K} \Cal N_a\sqcup\bigsqcup_{k\in\Bbb N}\bigsqcup_{b\in \Cal K}\Cal M_{b,k}
$$
of $\Bbb N$ into infinite subsets.
Recall that $\Cal K$ stands for the subset of $v$-periodic points in $K$.
For each $a\in\Cal K$, we denote  by $m_a$ the least positive period of $a$ under $v$.

Now we define a sequence $(C_n,h_n,z_n,\alpha_n)_{n=1}^\infty$ via an inductive procedure.
Suppose we have already constructed this sequence up to index $n$.
Consider now two cases.

Case {\bf [I]}. Let  $n+1\in \Cal N_a$ for some $a\in\Cal K$.
Now we set
$$
\gather
z_{n+1}:=2m_anh_n, \quad r_n:=n^3m_a,\\
C_{n+1}:=2h_n\cdot\{0,1,\dots, r_n-1\},\\
h_{n+1}:= 2r_nh_n,
\endgather
$$
Let $\alpha_{n+1}:C_{n+1}\to K$ be any map satisfying the following conditions
\roster
\item"(A1)" $\alpha_{n+1}(c+z_{n+1})=v\circ\alpha_{n+1}(c)$ for all $c\in C_{n+1}\cap (C_{n+1}-z_{n+1})$,
\item"(A2)" for each $0\le i< m_a$ there is a subset $C_{n+1,i}\subset C_{n+1}$ such that
$$
\gather
C_{n+1,i}-2h_n\subset C_{n+1}, \\  \alpha_{n+1}(c)=\alpha_{n+1}(c-2h_n)+v^i(a)\text{ for all }c\in C_{n+1,i}\text{  and}\\
\bigg|\frac{\#C_{n+1,i}}{\#C_{n+1}}-\frac 1{m_a}\bigg|<\frac 2{nm_a}.
\endgather
$$
\endroster
Case {\bf [II]}. Let  $n+1 \in \Cal M_{b,k}$ for some $b\in\Cal K$ and $k\in\Bbb N$.
 We set
$$
\gather
z_{n+1}:=m_{b}n(2h_n(k+1)+k), \quad r_n:=n^3(k+1)m_{b},\\
D_{n+1}:=2h_n\cdot\{0,1,\dots, nm_{b}-1\}\sqcup ((2h_n+1)\cdot\{1,2,\dots,nkm_{b}\}+2h_n(nm_{b}-1)),\\
C_{n+1}:=D_{n+1}+z_{n+1}\cdot\{0,1,\dots,n^2-1\},\\
h_{n+1}:= 2r_nh_n+kr_n/(k+1).
\endgather
$$
Let $\alpha_{n+1}:C_{n+1}\to K$ be any map satisfying the following conditions
\roster
\item"(A3)" $\alpha_{n+1}(c+z_{n+1})=v\circ\alpha_{n+1}(c)$ for all $c\in C_{n+1}\cap (C_{n+1}-z_{n+1})$,
\item"(A4)" for each  $0\le i< m_b$ there is a subset $C_{n+1,i}\subset C_{n+1}$ such that
$$
\gather
C_{n+1,i}-2h_n\subset C_{n+1}, \\  \alpha_{n+1}(c)=\alpha_{n+1}(c-2h_n)+v^i(b)\quad
\text{ for all }c\in C_{n+1,i}\text{  and}\\
\bigg|\frac{\#C_{n+1,i}}{\#C_{n+1}}-\frac 1{(k+1)m_{b}}\bigg|<\frac 2{nm_{b}}.
\endgather
$$
\item"(A5)"  there is a subset $C_{n+1,\vartriangle}\subset C_{n+1}$ such that
$$
\gather
C_{n+1,\vartriangle}-2h_n-1\subset C_{n+1}, \\  \alpha_{n+1}(c)=\alpha_{n+1}(c-2h_n-1)\quad
\text{ for all }c\in C_{n+1,\vartriangle}\text{  and}\\
\bigg|\frac{\#C_{n+1,\vartriangle}}{\#C_{n+1}}-\frac k{k+1}\bigg|<\frac 2n.
\endgather
$$
\endroster
Thus, $C_{n+1},h_{n+1}, z_{n+1},\alpha_{n+1}$ are completely defined.

We now let $F_n:=\{0,1,\dots,h_n-1\}$. Denote by $(X,\mu,T)$ the  $(C,F)$-transforma\-ti\-on associated with the sequence $(C_{n+1},F_n)_{n\ge 0}$.
Since
$$
\# F_{n+1}=h_{n+1}\ge 2r_nh_{n}=2\# C_{n+1}\# F_n,
$$
it follows that $\mu(X_{n+1})\ge 2\mu(X_n)$ for every $n$ and hence $\mu(X)=\infty$.

Let $\Cal R$ stand for the tail equivalence relation (or, equivalently, $T$-orbit equivalence relation) on $X$.
 Denote by $\alpha:\Cal R\to K$ the cocycle of $\Cal R$ associated with the sequence $(\alpha_n)_{n>0}$.

\head 4. Ergodicity and multiple recurrence of the transformations $T^{\times p}\times T_\alpha$
\endhead

\proclaim{Proposition 4.1}
The transformation $T^{\times p}\times T_\alpha$ is ergodic for each $p>0$.
\endproclaim
\demo{Proof}
We first show that $T^{\times p}$ is ergodic.
For simplicity, we will consider only the case when $p=2$.
(The general case is considered in a similar way.)
We note that for each  $n\in\Cal M_{b,1}-1$, there are subsets $C_{n+1}'$ and $C_{n+1}''$ in $C_{n+1}$
such that
\roster
\item"(a)"
if $c\in C_{n+1}'$ then $2h_n+c\in C_{n+1}$,
\item"(b)"
if $c\in C_{n+1}''$ then $2h_n+1+c\in C_{n+1}$,
\item"(c)"
$\# C_{n+1}'\ge 1/3\cdot\# C_{n+1}$ and  $\# C_{n+1}''\ge 1/3\cdot\# C_{n+1}$.
\endroster

Suppose we are given an arbitrary $n>0$ and $f,f',d,d'\in F_n$.
Assume for definiteness that $f\ge f'$ and $d\ge d'$.
Let $s:=\max(f-f',d-d')$.
We now find $k>n$ such that the intersection $\{n+1,\dots,k\}\cap\bigsqcup_{b\in\Cal K}(\Cal M_{b,1}-1)$
consists of $s$ points, say $l_1,\dots,l_s$.
We now set
$A:=f+\sum_{i=1}^k A_i$ and $B:=d+\sum_{i=1}^kB_i$, where
$$
A_i:=\cases
C_i &\text{if }i\notin\{l_1,\dots,l_s\}\\
C_i''&\text{if }i\in\{l_1,\dots,l_{f-f'}\}\\
C_i'&\text{if }i\in\{l_{f-f'+1},\dots,l_s\}
\endcases
\quad\text{and}\quad
B:=\cases
C_i &\text{if }i\notin\{l_1,\dots,l_s\}\\
C_i''&\text{if }i\in\{l_1,\dots,l_{d-d'}\}\\
C_i'&\text{if }i\in\{l_{d-d'+1},\dots,l_s\}
\endcases.
$$
It is easy to deduce from (a)--(c) that
\roster
\item"$(\circ)$"
$[A]_k\subset [f]_n$, $[B]_k\subset [d]_n$
\item"$(\circ)$"
$\mu([A]_k)\ge\frac 1{3^s}\mu([f]_n)$, $\mu([B]_k)\ge\frac 1{3^s}\mu([d]_n)$
\item"$(\circ)$"
$T^{2(h_{l_1-1}+\cdots+h_{l_s-1})}[A]_k\subset[f']_n$,
$T^{2(h_{l_1-1}+\cdots+h_{l_s-1})}[B]_k\subset[d']_n$.
\endroster
It remains to apply the following lemma.

\proclaim{Lemma 4.2} Fix $p>0$ and a map $\delta:\Bbb Z^p\to\Bbb R_+$ such that $\sum_{g\in\Bbb Z^p}\delta(g)<\frac 12$. If for each $n>0$ and $f_1,\dots,f_p,f_1',\dots,f_p'\in F_n$, there are a subset
$A\subset[f_1]_n\times\cdots\times[f_p]_n$ and $s\in\Bbb Z$ such that
$$
\mu^p(A)>\delta(f_1-f_1',\dots,f_p-f_p')
\mu^p([f_1]_n\times\cdots\times[f_p]_n)
$$
and $(T^{\times p})^sA\subset [f_1']_n\times\cdots\times[f_p']_n $
then  $T^{\times p}$ is ergodic.
\endproclaim

 This lemma is a particular case of \cite{Da1, Lemma~2.4} or \cite{DaS1, Lemma~5.2}.

Now we verify that the product $T^{\times p}\times T_\alpha$ is ergodic.
It is convenient to consider this transformation as a skew product
$(T^{\times (p+1)})_{1\otimes\alpha}$.
Fix $n>0$, $f_1,\dots,f_{p+1}\in F_n$ and $a\in\Cal K$.
We can find $k>n$ such that $k+1\in\Cal N_a$.
Then we  set
$$
A_i:=f_i+C_{n+1}+\cdots+ C_k+C_{k+1,0}
$$
and $A:=[A_1]_{k+1}\times\cdots\times[A_p]_{k+1}\subset[f_1]_n \times\cdots\times[f_p]_n$.
It follows from (A2) that
 $$
\gather
(T^{\times p})^{-2h_k}[A_1]_{k+1}\times\cdots\times[A_p]_{k+1}\subset
[f_1]_n \times\cdots\times[f_p]_n,\\
\frac{\mu^p(A)}{\mu^p([f_1]_n \times\cdots\times[f_p]_n)}>\bigg(\frac 1{2m_a}\bigg)^p\quad\text{and }\\
1\otimes\alpha(x, (T^{\times p})^{-2h_k}x)=\alpha(x_{p+1},T^{-2h_k}x_{p+1})=a
\endgather
$$
for all $x=(x_1,\dots,x_{p+1})\in A$.
Since $T^{\times(p+1)}$ is ergodic and $\Cal K$ is dense in $K$, we deduce from the standard ergodicity criterium for cocycles \cite{Sc} that the cocycle $1\otimes\alpha$ is ergodic, i.e.
$T^{\times p}\times T_\alpha$ is ergodic.
\qed
\enddemo

\remark{Remark \rom{4.3}}
We also note that $T^{\times p}$ is multiply recurrent for each $p>0$.
This follows from \cite{DaS1, Remark~2.4(i)}.
Now \cite{In} yields that $T^{\times p}\times T_\alpha$ is also multiply   recurrent.
 \endremark

\head 5. Proof of the main result (Theorem 0.1)
\endhead

Since the case $E=\{1\}$ is trivial (each infinite measure preserving rank-one transformation has a simple spectrum), we will assume from now on that $E\ne\{1\}$.
Then we fix an integer $m>0$ such that $m+1\in E$.
{\it Our purpose} is to prove that  $\Cal M(T^{\odot m}\times T_{\alpha,H})=E$, where the objects $T,\alpha,K,H$ were defined in the Section~3.

\proclaim{Lemma 5.1} Let  $a,b\in \Cal K$ and $k\in\Bbb N$.
Then for each $\chi\in\widehat K$,
\roster
\item"\rom{(i)}"
$U_{T,\chi}^{2h_n}\to l_\chi(a)\cdot I$
 as  $\Cal N_a-1\ni n\to\infty$ and
\item"\rom{(ii)}"
 $U_{T,\chi}^{2h_n}\to \frac{l_\chi(b)}{k+1}\cdot I
+ \frac{k}{k+1}U_{T,\chi}^*$ as  $\Cal M_{b,k}-1\ni n\to\infty$.
\endroster
\endproclaim

\demo{Proof} We show only (ii) since (i) is proved in a similar way but a bit simpler.
Let $n+1\in\Cal M_{b,k}$.

 Take any subset $A\subset F_n$.  We note that $[A]_n=[A+C_{n+1}]_{n+1}$.
Therefore it follows from (A4) that
for each $x\in T[F_n]_n$,
$$
\align
U_{T,\chi}^{2h_n}1_{[A]_n}(x)=&\sum_{i=0}^{m_{b}-1} \chi(\alpha(T^{2h_n}x,x)) 1_{[A+C_{n+1,i}]_{n+1}}(T^{2h_n}x)\\
&+\chi(\alpha(T^{2h_n}x,x))1_{[A+C_{n+1,
\vartriangle}]_{n+1}}(T^{2h_n}x)+\bar{o}(1)\\
=&\sum_{i=0}^{m_{b}-1} \chi(v^i(b)) 1_{[A+C_{n+1,i}-2h_n]_{n+1}}(x)\\
&+\chi(\alpha(T^{-1}x,x)) 1_{[A+C_{n+1,\vartriangle}-2h_n-1]_{n+1}}(T^{-1}x)+\bar{o}(1)\\
=&\sum_{i=0}^{m_{b}-1} \chi(v^i(b)) 1_{[A+C_{n+1,i}-2h_n]_{n+1}}(x)\\
&+U_{T,\chi}^* 1_{[A+C_{n+1,\vartriangle}-2h_n-1]_{n+1}}(x)+\bar{o}(1).\\
\endalign
$$
Therefore
$$
U_{T,\chi}^{2h_n}-\sum_{i=0}^{m_{b}-1} \chi(v^i(b))1_{[C_{n+1,i}-2h_n]_{n+1}}
-U_{T,\chi}^* 1_{[C_{n+1,\vartriangle}-2h_n-1]_{n+1}}\to 0
$$
as $\Cal M_{b,k}-1\ni n\to\infty$.
Here the functions $1_{[C_{n+1,i}-2h_n]_{n+1}}$ and  $1_{[C_{n+1,\vartriangle}-2h_n-1]_{n+1}}$ are considered as multiplication operators (orthogonal projectors) in $L^2(X,\mu)$.
It remains to  use the inequalities from (A4) and (A5) and a standard fact that for any sequence $C_{n}^\bullet\subset C_{n}$ such that $\# C_{n}^\bullet/\# C_{n}\to \delta$, we have
$$
1_{[C_{n}^\bullet]_{n}}\to \delta I\quad
\text{ as }n\to\infty.
$$
 \qed
\enddemo

\demo{Proof of  Theorem 0.1}
It follows from Proposition~4.1 that the transformation $T^{\odot m}\times T_{\alpha,H}$ is ergodic.
There is a natural decomposition
of   $U_{T^{\odot m}\times T_{\alpha,H}}$ into
 an orthogonal direct sum
$$
U_{T^{\odot m}\times T_{\alpha,H}}=\bigoplus_{\chi\in \widehat {K/H}}(U_T^{\odot m}\otimes U_{T,\chi}).
$$
Let us show the following claims:
\roster
\item"(i)"
$U_T^{\odot m}\otimes U_T$ has a homogeneous spectrum of multiplicity $m+1$,
\item"(ii)"
$U_T^{\odot m}\otimes U_{T,\chi}$ has a simple spectrum if $\widehat{K/H}\ni\chi\ne 0$,
\item"(iii)"
$U_T^{\odot m}\otimes U_{T,\chi}$ and $U_T^{\odot m}\otimes U_{T,\xi}$ are unitarily equivalent if $\chi$ and $\xi$ belong to the same $\widehat v$-orbit,
\item"(iv)" the measures of maximal spectral type of $U_T^{\odot m}\otimes U_{T,\chi}$ and $U_{T}^{\odot m}\otimes U_{T,\xi}$ are mutually singular if $\chi$ and $\xi$ do not belong to the same $\widehat v$-orbit.
\endroster

By Lemma~5.1(ii),
$$
U_T^{2h_n}\to\frac 1{k+1}I+\frac k{k+1}U_T^*\tag5-1
$$
as $\Cal M_{b,k}-1\ni n\to\infty$ for each  $b\in\Cal K$ and $k\in\Bbb N$.
If follows from the proof of Lemma~2.2 that the unitary operator $U_T^{\odot (m+1)}$ has a simple spectrum.
Now the claim (i) follows from Corollary~1.2.

Since $T$ is of rank one and the map $[f]\ni x\mapsto \alpha(Tx,x)\in K$ is constant for each $f\in F_n\setminus\{h_n-1\}$, $n\in\Bbb N$, it follows that the operator $U_{T,\chi}$ has a simple spectrum for each $\chi\in\widehat{K/H}$.
If $\widehat{K/H}\ni\chi\ne 0$ then by
 Lemma~2.1(ii), there is  $b\in\Cal K$ such that the numbers $l_\chi(b)\ne 1$.
Lemma~5.1(ii) yields
$$
U_{T,\chi}^{2h_n}\to \frac {l_\chi(b)}{k+1}\cdot I
+ \frac k{k+1}U_{T,\chi}^*\tag5-2
$$
as $\Cal M_{b,k}-1\ni n\to\infty$, $k=1,\dots,m$.
Applying Lemma~2.2 with \thetag{5-1} and \thetag{5-2} we obtain the claim (ii).

Since
$$
\sum_{n>0}
\frac{\#(C_n\triangle (C_n-z_n))}{\# C_n}=\sum_{n>0}\frac 2{n^2},
$$
it follows from Lemma 3.1 that a transformation $S_{ \bar z}$ of $(X,\mu)$ is well defined by the formula \thetag{3-1} and $S_{\bar z}\in C(T)$.
 It follows from (A1) and (A3) that \thetag{3-2} is satisfied. Hence by Lemma~3.1,
$$
\alpha\circ S_{\bar z}\quad\text{is cohomologous to }v\circ\alpha
$$
and (iii) follows from Remark~2.3.

Let characters $\chi,\xi\in\widehat{K}$ do not belong to the same $\widehat v$-orbit.
By Lemma~2.1(ii), there is $a\in\Cal K$ with $l_\chi(a)\ne l_{\xi}(a)$.
We now deduce from Lemma~5.1(i) that
$$
(U_T^{\odot m}\otimes U_{T,\chi})^{2h_n}\to l_\chi(a) I\quad{\text{ and \ }}
(U_T^{\odot m}\otimes U_{T,\xi})^{2h_n}\to l_\xi(a) I
$$
as  $\Cal N_{a}-1\ni n\to\infty$.
This yields (iv).

Now (i)--(iv) imply
$$
\Cal M(T^{\odot m}\times T_{\alpha,H})=\{m+1\}\cup E=E.
$$

\qed
\enddemo

\enddocument